\documentclass[11pt]{article}
\usepackage{amsmath,amsfonts,amssymb}
\usepackage{bbm}
\usepackage{cases}
\usepackage[hidelinks]{hyperref}
\usepackage{xcolor}

\def\<{{\langle}}
\def\>{{\rangle}}
\allowdisplaybreaks[4]

\newtheorem{theorem}{Theorem}[section]
\newtheorem{lemma}[theorem]{Lemma}
\newtheorem{corollary}[theorem]{Corollary}
\newtheorem{proposition}[theorem]{Proposition}

\newtheorem{remark}[theorem]{Remark}

\newtheorem{definition}[theorem]{Definition}
\newtheorem{assumption}[theorem]{Assumption}

\begin{document}
	
	\title{\vspace*{-1.5cm}
		 Stochastic heat equations on moving domains}
	
	\author{Tianyi Pan$^{1}$, Wei Wang$^{1}$, Jianliang Zhai$^{1}$,
		Tusheng Zhang$^{2}$}
	\footnotetext[1]{\, School of Mathematics, University of Science and Technology of China, Hefei, China. Email: pty0512@mail.ustc.edu.cn (Tianyi Pan), ww12358@mail.ustc.edu.cn (Wei Wang), zhaijl@ustc.edu.cn (Jianliang Zhai).}
	\footnotetext[2]{\, School of Mathematics, University of Manchester, Oxford Road, Manchester M13 9PL, England, U.K. Email: tusheng.zhang@manchester.ac.uk}
	\date{}
	\maketitle
	
	\begin{abstract}
		In this paper, we establish the well-posedness of stochastic heat equations on moving domains, which amounts to a study of infinite dimensional interacting systems. The main difficulty is to deal with the problems caused by the time-varying state spaces and the interaction of the particle systems. The interaction still occurs even in the case of additive noise. This is in contrast to stochastic heat equations in a fixed domain.
	\end{abstract}
	\noindent
	{\bf Keywords and Phrases:} Stochastic partial differential equations,  stochastic heat equations, moving domain, space-time noise.
	
	\medskip
	
	\noindent
	{\bf AMS Subject Classification:} Primary 60H15;  Secondary 35R60.
	\section{Introduction}
	$\qquad$In this paper, we are concerned with the well-posedness of stochastic heat equations driven by multiplicative noise on a moving domain, which is formally given as follows,
	\begin{numcases}{}\label{SHE}
		\frac{\partial u(t,x)}{\partial t}=\Delta_t u(t,x)+\sigma\big(t,u(t,x)\big)\frac{\partial^2W(t,x)}{\partial t\partial x}\quad \text{ for }(t,x)\in \mathcal{O}_T,\\
		u(t,x)=0\quad\text{ for }x\in \partial I_t, \nonumber\\
		u(0,x)=u_0(x)\quad\text{ for } x\in I_0\nonumber.
	\end{numcases}
Here $\mathcal{O}_T=\mathop{\cup}\limits_{0\leq t\leq T} \{t\}\times (0,a_t)$ is a non-cylindrical time-space domain, $ a:[0,T]\rightarrow (0,+\infty)$  is a continuously differentiable function, $\ \Delta_t$ stands for the Dirichlet Laplacian operator on $I_t:=(0,a_t)$ with boundary $\partial I_t=\{0,a_t\}$, and $\frac{\partial^2W(t,x)}{\partial t\partial x}$ stands for a space-time noise specified later.
	\vskip 0.3cm
	The classical heat equation, i.e. equation (\ref{SHE}) with $\sigma\equiv0$ and $I_t\equiv I_0$, describes the evolution of the heat flow on a fixed domain. The solution $u(t,x)$ is the heat density at position $x\in I_0$ and time $t\in [0,T]$ with initial density $u_0(x)$. In the stochastic setting, the noise term represents  a random internal/external heat source.  The equation (\ref{SHE}) describes the time evolution of the heat density in  a domain moving with time and  with a time-dependent random heat source.\par
    To solve partial differential equations(PDEs) on time-varying domains, there are mainly two classes of methods in the literature. One is called the ``diffeomorphism method". With this method the original PDEs on the moving domain is pushed into  PDEs on a flat domain via a family of diffeomorphisms. Although the new PDEs are on a fixed domain, they have  much  higher nonlinear terms than the original PDEs. We refer the readers to \cite{MT,IW,HH, FKW} and references therein for details. The other is called the ``penalty method". The PDEs are solved via a penalty term. We also refer the readers to \cite{HN,KSR,BG} for details.  \par
    For the stochastically-forced equations on time-varying domains, the problem becomes much more delicate because of the singularities introduced by the noise term.  There are very few results on the well-posedness of stochastic partial differential equations on moving domains.  We refer the reader to a recent paper \cite{WZZ} by some of the authors of this paper, in which the authors established  well-posedness for stochastic 2D Navier-Stokes equations on a time-dependent domain driven by an additive noise. In that paper, the ``diffeomorphism method" has been adopted.  However, using this method, it seems very difficult to built the $\rm{It\hat{o}}$ formula for the solutions, and  to  solve  the case of multiplicative noises. Similar problems arise when one applies the ``penalty method". The methods mentioned above are suitable to solve PDEs on time-varying domains, but it seems less ineffective when one uses them to deal with the case of the stochastically-forced equations.

The purpose of this paper is to establish the well-posedness of stochastic heat equations  driven by multiplicative noise on one dimensional moving domains. We use the Galerkin method, which is a commonly used method to obtain the well-posedness of (stochastic) PDEs on fixed domains, but this is the first paper to use this method to deal with the case of stochastic equations on moving domains, even the case of deterministic equations on moving domains. Compared with the proof of the corresponding results for the case of fixed
domains, new essential difficulties appear. For example, the eigenbasis of time-dependent Laplacian $\Delta_t$ is dependent on $t$. We now describe the approach of this paper in more details. Using the evolving eigenbasis of time-dependent Laplacian, it turns out that the problem becomes a study of the well-posedness of an infinite interacting particle system. We then use  finite dimensional approximations. Through proving the tightness of the laws of the approximating solutions, we first obtain the existence of a probabilistic weak solution to equation (\ref{SHE}) and then the well-posedness of strong solutions by combining the pathwise uniqueness and Yamada-Watanable theorem. The main difficulty is to deal with the problem caused by the time-varying state space and the interaction of the particle systems. We notice that the corresponding system is still interacting even in the case of additive noise. This is in contrast to stochastic heat equations on a fixed domain. Compared with the ``diffeomorphism method" and the ``penalty method" mentioned above, an important novelty of our approach is that we are able to get  an $\rm{It\hat{o}}$ formula/energy identity for the solutions and to handle
the case of multiplicative noise.
Our method is robust enough to handle stochastic PDEs with highly nonlinear terms, e.g. stochastic porous medium equation,
stochastic $p$-Laplace equation, stochastic Burgers type equations, stochastic 2D Navier-Stokes equations, which will be the subject of our future work.

%
%

    \vskip 0.3cm

    The rest of the paper is organized as follows. In Section 2, we provide the framework and state the main result. In Section 3, we study the approximating solutions. Among other things we show that the approximating solutions is a Cauchy sequence on a new probability space.  Section 4 is devoted to the proof of the main result.\par
    \section{Assumptions and main results}
    \ \ \ \ In this section, we present the notations, the precise assumptions and the main result. Recall that we assume that the domain  $I_t$ deforms in time $t$ in a continuously differentiable way, more precisely,
    \begin{assumption}\label{assump1}
    	There exists a $\mathcal{C}^1$ function $a:[0,T]\rightarrow(0,\infty)$ such that $I_t=(0,a_t)$, $\forall \ t\in[0,T]$.
    \end{assumption}
    In the following, the derivative of $a_t$ is denoted by $a'_t$.
\begin{remark}\label{Rem 1}
	Under assumption \ref{assump1},  there exist strictly positive constants $\delta_0$ and $L$ such that $\delta_0\leq a_t\leq L$ and $|a'_t|\leq L$, $\forall \ t\in[0,T]$.
\end{remark}\par
  We let $\mathbb{L}^2(I_t)$ be the space of all square-integrable functions on $I_t$ and $\mathbb{H}^1_0(I_t)$  the closure of the space of all smooth functions compactly supported in $I_t$ under the norm
    $$\|u\|_{\mathbb{H}_0^1(I_t)}=\|\partial_xu\|_{\mathbb{L}^2(I_t)}.$$
	For simplicity, $\forall\ t\geq0$, we write $\|\cdot\|_{\mathbb{L}^2(I_t)}$  as $|\cdot|_t$, $(\cdot,\cdot)_{\mathbb{L}^2(I_t)}$  as $(\cdot,\cdot)_t$, $\|\cdot\|_{\mathbb{H}_0^1(I_t)}$  as $\|\cdot\|_t$, and  $(\cdot,\cdot)_{\mathbb{H}_0^1(I_t)}$  as $\langle\cdot,\cdot\rangle_t$. Sometimes we will also write $\mathbb{L}^2(I_t)$ and $\mathbb{H}^1_0(I_t)$ as $\mathbb{H}_t$ and $\mathbb{V}_t$ respectively.\par

	Now we introduce the function spaces that the solution belongs to.
		For $L>0$ appearing in Remark \ref{Rem 1}, the function space $\mathbb{L}^2(I_t)$ can be isometrically embedded into $\mathbb{H}:=\mathbb{L}^2\left((0,L)\right)$ by setting $f(t)\equiv0$ for $t\in(a_t,L)$ for any $f\in \mathbb{L}^2(I_t)$.
\vskip 0.3cm
 Define
		\begin{align*}
			\mathbb{X}:=\Big\{f\in C\left([0,T];\mathbb{H}\right),f(t)\in \mathbb{L}^2(I_t),\ \forall\ t\in[0,T]\Big\}
		\end{align*}
		equipped with the norm
		\begin{align*}
			\|f\|_\mathbb{X}=\sup\limits_{0\leq t\leq T}{|f(t)|_t}.
		\end{align*}
		
\noindent It can be easily verified that $\mathbb{X}$ is a Banach space, isometrically embedded into $C\left([0,T];\mathbb{H}\right)$.

		Similarly, we can isometrically embed $\mathbb{H}_0^1(I_t)$ into $\mathbb{V}:=\mathbb{H}_0^1\left((0,L)\right)$ by setting the functions in $\mathbb{H}_0^1(I_t)$ to be zero outside the interval $I_t$.  This extension is well-defined since the functions in $\mathbb{H}_0^1(I_t)$ vanish at the boundary. We define the space:
		\begin{align*}
			\mathbb{Y}:=\Big\{f\in \mathbb{L}^2\left([0,T];\mathbb{H}_0^1((0,L))\right),f(t)\in \mathbb{H}_0^1(I_t)\ a.e.\ t\in[0,T]\Big\}
		\end{align*}
		equipped with the norm
		\begin{align*}
			\|f\|^2_\mathbb{Y}=\int_{0}^{T}\|f(t)\|^2_t\mathrm{d}t.
		\end{align*}
		It can also be verified that $\mathbb{Y}$ is a Hilbert space, isometrically embedded into $\mathbb{L}^2\left([0,T];\mathbb{V}\right)$.
We will  write $\|\cdot\|_\mathbb{H}$, $(\cdot,\cdot)_\mathbb{H}$, $\|\cdot\|_\mathbb{V}$, $(\cdot\ ,\cdot)_\mathbb{V}$  as $|\cdot|$, $(\cdot\ ,\cdot)$,  $\|\cdot\|$ $\langle\cdot,\cdot\rangle$, respectively.\par

Analogously, we introduce the space:
		\begin{align*}
			\mathbb{Z}:=\Big\{f\in \mathbb{L}^2\left([0,T];\mathbb{H}\right),f(t)\in \mathbb{L}^2(I_t)\ a.e.\ t\in[0,T]\Big\}
		\end{align*}
		equipped with the norm
		\begin{align*}
			\|f\|^2_\mathbb{Z}=\int_{0}^{T}|f(t)|^2_t\mathrm{d}t.
		\end{align*}

    For the Dirichlet Laplacian $\Delta_t$ on $I_t$, there exists an orthonormal eigensystem $\big\{\lambda_k(t),e_k(t)\big\}_{k\geq1}$ on $\mathbb{H}_t$ given by
\begin{equation}\label{eigen}
	\lambda_k(t)=-\left(\frac{k\pi}{a_t}\right)^2; \quad e_k(t,x)=\sqrt{\frac{2}{a_t}}\sin \frac{k\pi x}{a_t},\ x\in I_t.
\end{equation}
 $\{e_k(t)\}_{k\in\mathbb{N}}$ is an orthonormal basis of $\mathbb{H}_t$ and  also an orthogonal basis in $\mathbb{V}_t$.
\vskip 0.3cm
		Take a sequence of independent standard real-valued Brownian motions $\{B^k\}_{k\geq 1}$ on a filtered probability space $(\Omega,\mathcal{F},\{\mathcal{F}_t\}_{t\geq 0},P)$ satisfying the usual conditions, and an orthonormal basis $\{f_k\}_{k\geq1}$ of $\mathbb{L}^2\big((0,L)\big)$, then we introduce a cylindrical Brownian motion $W$ on $\mathbb{H}$ given by $$W_t:=\sum_{k=1}^{\infty}f_kB^k_t, \ t\geq 0.$$\par
    Now we give the assumptions on the mapping $\sigma$ appearing in the stochastic heat equation (\ref{SHE}).
    Let $\sigma:[0,T]\times \mathbb{H}\rightarrow\mathbb{L}_2(\mathbb{H},\mathbb{H})$ be a Borel measurable mapping, where $\mathbb{L}_2(\mathbb{U}_1,\mathbb{U}_2)$ denotes  the space of Hilbert-Schimdit operators  from a Hilbert space $\mathbb{U}_1$ to another Hilbert space $\mathbb{U}_2$ equipped
with the usual Hilbert-Schmidt norm  $\|\cdot\|_{\text{HS}(\mathbb{U}_1,\mathbb{U}_2)}$. When there is no danger of causing ambiguity,
we write $\|\cdot\|_{\text{HS}}=\|\cdot\|_{\text{HS}(\mathbb{U}_1,\mathbb{U}_2)}$.
   \begin{assumption}\label{assump2}
 We  assume that $\sigma$ is Lipschitz continuous and measurable, i.e., there is a positive constant $K>0$ such that, for any $s\in[0,T]$ and
 $u,v,h\in \mathbb{H}_s$,
    \begin{itemize}
    	\item[\rm (i)]\ $\|\sigma(s,u)-\sigma(s,v)\|_{\rm{HS}}\leq K\|u-v\|_\mathbb{H}$;
    	\item [\rm  (ii)]\ $\|\sigma(s,u)\|_{\rm{HS}}\leq K(\|u\|_\mathbb{H}+1)$;
    	\item [\rm  (iii)]\ $\sigma(s,h)\in\mathbb{L}_2(\mathbb{H},\mathbb{H}_s)$.
    \end{itemize}
\end{assumption}
The stochastic integral $\sigma\big(t,u(t)\big)\text{d}W_t$ is understood as
$$\sigma\big(t,u(t)\big)\text{d}W_t=\sum_{i=1}^{\infty}\sigma\big(t,u(t)\big)f_k\text{d} B_t^k.$$
 	Let us denote the time-derivative  $\partial_s\varphi(s,x)$ by $\varphi'(s,x)$. Recall the definition of $\mathcal{O}_T$  in Section 1, we have the following definition of the solution to (\ref{SHE}).
	\begin{definition}\label{sol1}
		$u\in \mathbb{X}\ a.s.$ is a solution of equation (\ref{SHE}) if it is predictable with respect to $\{\mathcal{F}_t\}_{t\geq0}$ as an $\mathbb{H}$-valued random process and, for $\forall\ t\in [0,T]$ and $\varphi\in \mathcal{C}_0^\infty (\bar{\mathcal{O}}_T)=\big\{g\in \mathcal{C}^\infty(\mathcal{\bar{O}}_T):\text{$g(t,0)=g(t,a_t)=0, \forall t\in[0,T]$}$\big\},
		\begin{align}\label{sol}
			&\int_{0}^{a_t} u(t,x)\varphi(t,x)\mathrm{d}x -\int_{0}^{a_0} u_0(x)\varphi(0,x) \mathrm{d}x-\int_{0}^{t}\int_{0}^{a_s} u(s,x)\varphi'(s,x) \mathrm{d}x\mathrm{d}s\nonumber\\
			=&\int_{0}^{t}\int_{0}^{a_s}  u(s,x)\Delta_s\varphi(s,x)\mathrm{d}x\mathrm{d}s+\int_0^t\Big(\varphi (s),\sigma\big(s,u(s)\big)\mathrm{d}W_s\Big)_{\mathbb{H}}  \ a.s.
		\end{align}
	\end{definition}
	
	Now we are ready to state the main result.
	\begin{theorem}\label{WP}
		For a deterministic function $u_0\in\mathbb{L}^2\left(I_0\right)$, there exists a unique solution $u\in \mathbb{L}^2(\Omega;\mathbb{X}\cap \mathbb{Y})$ to equation (\ref{SHE}).
	\end{theorem}
\par
Throughout the paper, $ C$ denotes a generic  positive constant whose value may change from line to line. All other constants will be denoted by $C_1,C_2,\cdots$. They are all positive but their values are not important. The dependence of constants on parameters if needed will be indicated, e.g., $C_T$.\par
	\section{Approximating solutions}
In this section, we will study the approximating equations (\ref{u^n}) and provide a number of uniform estimates for the approximating solutions. Moreover, we will show that the approximating solutions form a Cauchy sequence on a new probability space. This plays an important role in the proof of the main result.
\vskip 0.3cm
	To formulate the appropriate approximating  equations (\ref{u^n}), we start with some formal calculations. Expand the solution $u$ (if it exists) of equation (\ref{SHE}) with respect to the eigenbasis $e_k(s,x),\ k\in\mathbb{N}$ (see (\ref{eigen})) to get
$$ u(s,x)=\sum_{k=1}^{\infty}\big(u(s),e_k(s)\big)_s e_k(s,x)=\sum_{k=1}^{\infty}A_k(s)e_k(s,x),$$
here
$$A_k(s):=\big( u(s),e_k(s)\big)_s=\int_0^{a_s}u(s,x)e_k(s,x) \mathrm{d}x.$$
	Taking the eigenbasis $e_k(t,x), k\in\mathbb{N}$ as the test functions in (\ref{sol}),
	it follows that
	\begin{align*}
		&\ \ \ \ A_k(t)-(u_0,e_k(0))_0-\int_{0}^{t}\int_{0}^{a_s} u(s,x) e_k'(s,x) \mathrm{d}x\mathrm{d}s
		\\&= \int_{0}^{t}\lambda_k(s)A_k(s)\mathrm{d}s+\int_{0}^{t}\Big(e_k(s),\sigma\big(s,u(s)\big)\mathrm{d}W_s\Big)_\mathbb{H}.
	\end{align*}
	A formal calculation yields that
	\begin{align}\label{akeq0}
		&\ \ \ \ A_k(t)-(u_0,e_k(0))_0-\sum_{j=1}^{\infty}\int_{0}^{t}A_j(s) b_{jk}(s) \mathrm{d}s
		\\&=\int_{0}^{t}\lambda_k(s)A_k(s)\mathrm{d}s+\int_{0}^{t}\Big(e_k(s),\sigma\big(s,\sum_{j=1}^{\infty}A_j(s)e_j(s)\big)\mathrm{d}W_s\Big)_\mathbb{H}\nonumber,
	 \end{align}
	where
	\begin{align}
		b_{jk}(s)=\int_0^{a_s} e_j(s,x)e_k^{\prime}(s,x)\mathrm{d}x=
		\begin{cases}
			(-1)^{j+k}\cdot\frac{a_s^{\prime}}{a_s}\cdot \frac{2jk}{j^2-k^2}, \qquad &j\neq k,\\
			0, \qquad &j=k.
		\end{cases}
	\end{align}
This suggests to  consider the following interacting systems of  stochastic  differential equations: for $k=1,\cdots,n$,
	\begin{align}\label{approxeq}
	&\ \ \ \ A^n_k(t)-(u_0,e_k(0))_0-\sum_{j=1}^{n}\int_{0}^{t}A^n_j(s) b_{jk}(s) \mathrm{d}s
	\\&=\int_{0}^{t}\lambda_k(s)A^n_k(s)\mathrm{d}s+\int_{0}^{t}\Big(e_k(s),\sigma\big(s,\sum_{j=1}^{n}A^n_j(s)e_j(s)\big)\mathrm{d}W_s\Big)_\mathbb{H}.\nonumber
\end{align}
	Note that the SDE (\ref{approxeq}) has a unique $\mathbb{R}^n$-valued continuous  solution $\big(A_1^n(t),\\
A_2^n(t),...,A_n^n(t)\big)$ since all the coefficients are Lipschitz continuous due to assumptions \ref{assump1} and \ref{assump2} (i). Set
	\begin{align}\label{u^n}
		u^n(t,x):=\sum\limits_{k=1}^{n}A_k^n(t)e_k(t,x).
	\end{align}
We begin with some moment estimates of $\{u^n\}_{n\geq1}$ in the space $\mathbb{X}\cap \mathbb{Y}$.
	\begin{proposition}\label{estimate}
		There exists a positive constant $C_1>0$ such that
		$$\sup_{n\in \mathbb{N}}E\Big[\sup\limits_{0\leq t\leq T}\big|u^n(t)\big|^2_t+\int_{0}^{T}\big\|u^n(t)\big\|_t^2\mathrm{d}t\Big]\leq C_1.$$
	\end{proposition}
	\textbf{Proof.}
	By ${\rm It\hat{o}}$'s formula, for  $k\leq n$,
	\begin{align*}
		|A_k^n(t)|^2=&\ |(u_0,e_k(0))_0|^2+2\sum_{j=1}^{n} \int_0^t A_j^n(s)b_{jk}(s) A_k^n(s)\mathrm{d}s+2\int_0^t \lambda_k(s)|A_k^n(s)|^2 \mathrm{d}s\\&+2\int_0^t\Big(A_k^n(s)e_k(s),\sigma\big(s,\sum_{l=1}^{n}A^n_l(s)e_l(s)\big)\mathrm{d}W_s\Big)\\&+\int_{0}^{t}\sum_{j=1}^{\infty}\Big|\sigma^k_j\big(s,\sum_{l=1}^{n}A^n_l(s)e_l(s)\big)\Big|^2\mathrm{d}s,
	\end{align*} where  $\sigma_j^k(s,u):=\big(e_k(s),\sigma(s,u)f_j\big)_s$. It follows that
	\begin{align}\label{un2}
		|u^n(t)|^2_{t}=&\sum_{k=1}^{n} |A_k^n(t)|^2 =\sum_{k=1}^{n}|(u_0,e_k(0))_0|^2+2\sum_{k=1}^{n}\sum_{j=1}^{n} \int_0^t A_j^n (s)b_{jk}(s) A_k^n(s)\mathrm{d}s\nonumber\\
		\ &+2\sum_{k=1}^{n} \int_0^t \lambda_k(s)|A_k^n(s)|^2 \mathrm{d}s+2\int_0^t\Big(u^n(s),\sigma\big(s,u^n(s)\big)\mathrm{d}W_s\Big)\nonumber\\&+\int_{0}^{t}\sum_{k=1}^{n}\sum_{j=1}^{\infty}\big|\sigma^k_j\big(s,u^n(s)\big)\big|^2\mathrm{d}s.
	\end{align}
	The second term on the right hand side is 0 since ${b_{jk}(s)}$ is skew-symmetric with respect to  $(j,k)$.
	Then by the BDG inequality and Young inequality, we have
	\begin{eqnarray*}
		&&E\left\{ \sup_{0\leq s \leq t} \big|u^n(s)\big|_s^2+2\int_0^t \big\|u^n(s)\big\|^2_s \mathrm{d}s\right\}\\
		&\leq &\ \sum_{k=1}^{n} |(u_0,e_k(0))_0|^2+CE\left\{  \left(\int_0^t \big|u^n(s)\big|_s^2 \big\|\sigma\big(s,u^n(s)\big)\big\|_{\text{HS}}^2\mathrm{d}s\right)^{1/2}\right\}\\
		&&+E\Big[\int_{0}^{t}\big\|\sigma\big(s,u^n(s)\big)\big\|_{\text{HS}}^2\mathrm{d}s\Big] \\
		&\leq &\ \sum_{k=1}^{n} |(u_0,e_k(0))_0|^2+C_K\Big(E \big[\sup_{0\leq s \leq t}\big|u^n(t)\big|_s^2\big]\Big)^{\frac{1}{2}}\left(E\int_{0}^{t}\big(1+|u^n(s)|^2_s\big)\mathrm{d}s\right)^{\frac{1}{2}}\\
		&&+2K^2t+2K^2E\Big[\int_{0}^{t}|u^n(s)|^2_s\mathrm{d}s\Big]\\
		&\leq& \ |u_0|_0^2+C_{T,K}+\frac{1}{2}E\big[ \sup_{0\leq s \leq t}|u^n(s)|_s^2\big]+C_K\int_{0}^{t}E\big[|u^n(s)|^2_s\big]\mathrm{d}s.
	\end{eqnarray*}
	By the Gronwall inequality, we obtain
	\begin{align*}
	E \Big[\sup_{0\leq t \leq T} \big|u^n(t)\big|_{t}^2+2\int_0^T \|u^n(s)\|^2_s \mathrm{d}s\Big]\leq \Big(|u(0)|_0^2+C_{T,K}\Big)C_{T,K}.
	\end{align*}
	Since the right hand side is independent of $n$, the proof is complete. $\blacksquare$
\vskip 0.3cm
We also have the uniform estimate of higher order moments.
\begin{corollary}\label{moest3}
	For all $p>1$ there exists a positive constant $C_{p,K,T}$ such that $$\sup_{n\in \mathbb{N}} E\Big\{\sup_{0\leq t \leq T} \big|u^n(t)\big|_t^{2p}\Big\}\leq C_{p,K,T},$$
	 $$\sup_{n\in \mathbb{N}} E\Big\{\Big(\int_{0}^{T} \|u^n(t)\|_t^{2}\mathrm{d}t\Big)^p\Big\}\leq C_{p,K,T}.$$
\end{corollary}
\textbf{Proof.} By (\ref{un2}), we find that
\begin{eqnarray*}
	&& \sup_{0\leq s \leq t}|u^n(s)|^{2p}_{s}+\left(\int_{0}^{t}\|u^n(s)\|_s^2\mathrm{d}s\right)^p\\
	& \leq& C_p\Bigg(|u_0|_0^{2p}+\sup_{0\leq s \leq t}\Big|\int_0^s\Big(u^n(r),\sigma\big(r,u^n(r)\big)\mathrm{d}W_r\Big)\Big|^p+\Big(\int_{0}^{t}\sum_{k=1}^{n}\sum_{j=1}^{\infty}\Big|\sigma^k_j\big(s,u^n(s)\big)\Big|^2\mathrm{d}s\Big)^p\Bigg).
\end{eqnarray*}
 We define the truncating stopping time $$\tau_M^n:=\inf\big\{t>0:|u^n(t)|_t>M\big\}\wedge T.$$
 Replacing $t$ by $t\wedge\tau_M^n$ in the above inequality,  applying BDG's inequality and using  Young's inequality and Assumption \ref{assump2} (ii), we find
 \begin{align}\label{moest2}
 	&\ E\Big[\sup_{0\leq s \leq t\wedge\tau_M^n}|u^n(s)|^{2p}_{s}+\Big(\int_{0}^{t\wedge\tau_M^n}\big\|u^n(s)\big\|_s^2\mathrm{d}s\Big)^p\Big]\nonumber\\
	\leq&\  C_p\left(\!|u_0|_0^{2p}+\!E\Big\{\Big|\int_0^{t\wedge\tau_M^n}\big|u^n(s)\big|^{2}_{s}\big\|\sigma\big(s,u^n(s)\big)\big\|^2_{\text{HS}}\mathrm{d}s\Big|^\frac{p}{2}\!\!+\!\Big(\!\int_{0}^{t\wedge\tau_M^n}\big\|\sigma\big(s,u^n(s)\big)\big\|_{\text{HS}}^2\mathrm{d}s\Big)^p\Big\}\!\right)\nonumber\\
	\leq&\ C_p|u_0|_0^{2p}+\frac{1}{2}E\Big[\sup_{0\leq s \leq t\wedge\tau_M^n}|u^n(s)|^{2p}_{s}\Big]+C_{p,T,K}+C_{p,T,K}\int_{0}^{t}E\big[|u^n(s\wedge\tau_M^n)|_{s\wedge\tau_M^n}^{2p}\big]\mathrm{d}s.
 \end{align}
Then from the Gronwall inequality, we can find a constant $C_{p,K,T}>0$ independent of $n,M$ such that
$$E\Big[\sup_{0\leq s \leq T\wedge\tau_M^n}|u^n(s)|^{2p}_{s}+\Big(\int_{0}^{T\wedge\tau_M^n}\big\|u^n(t)\big\|_t^2\mathrm{d}t\Big)^p\Big]\leq C_{p,K,T}.$$
By letting $M\rightarrow\infty$ we obtain the desired estimates.
$\hfill\blacksquare$\\
\vskip 0.2cm
As a consequence, we obtain the following
\begin{corollary}\label{probest}
	For any $\varepsilon>0$, there exists $M>0$ such that $$\sup_{n\in\mathbb{N}}P\Big(\sup_{t\in[0,T]}\big|u^n(t)\big|_t>M\Big)\leq\varepsilon,$$
	$$\sup_{n\in\mathbb{N}}P\Big(\int_{0}^{T}\big\|u^n(t)\big\|^2_t\mathrm{d}t>M\Big)\leq\varepsilon.$$
\end{corollary}

    To prove the convergence of ${u^n}, n\geq 1$, we first need a compactness result. For this, we introduce a lemma about compact embedding proved in \cite{WZZ}.
\vskip 0.3cm
    Let $J_i$ denote a family of equicontinuous real-valued functions on $[0,T]$. For a positive constant $N$, set $K_{N,J}:=\mathop{\cap}^{\infty}_{i=1}K_{N,J_i}$, where
    \begin{align*}
    	K_{N,J_i}&=\Bigg\{g\in \mathbb{X}\cap \mathbb{Y}:\!\sup_{0\leq t \leq T}|g(t)|_t\leq N,\!\int_{0}^{T}\|g(t)\|^2_t\mathrm{d}t\leq N,\\&\ \ \ \ g_i=\big\{g_i(t):=\big(g(t),e_i(t)\big)_t,\ t\in[0,T]\big\}\in{J_i}\Bigg\}.
    \end{align*}
    The following lemma was proved in \cite{WZZ}.
	\begin{lemma}\label{tight}
		$K_{N,J}$ is precompact in $\mathbb{Z}$.
	\end{lemma}
	 Next we will establish the tightness of $\{\mathcal{L}(u^n)\}_{n\geq1}$, the family of distributions of $\{u^n\}_{n\in\mathbb{N}}$.
	
    \begin{proposition}
    	$\{\mathcal{L}(u^m)\}_{m\geq1}$ is tight in $\mathbb{Z}$.
    \end{proposition}
	 \textbf{Proof.}
	 Consider $k,m\in\mathbb{N}$ with $k\leq m$. Note that for any $0\leq s\leq t\leq T$,
	 \begin{eqnarray*}
	 &&A_k^m(t)-A_k^m(s)\\
	 &=&\sum_{j=1}^{m}\int_{s}^{t}A_j^m(r)b_{jk}(r)\mathrm{d}r+\int_{s}^{t}\lambda_k(r)A_k^m(r)\mathrm{d}r+\int_{s}^{t}\Big(e_k(r),\sigma\big(r,u^m(r)\big)\mathrm{d}W_r\Big)\\
       &=:&I_k^m(t,s) + II_k^m(t,s) + III_k^m(t,s).
	 \end{eqnarray*}
 We will estimate each of the terms on the right. For $p\in(2,\infty)$, by the BDG inequality we have
\begin{align*}
	E\big[|III_k^m(t,s)|^{\frac{8p}{p-2}}\big]&\leq C_pE\Big[
	\Big(\int_{s}^{t}\big\|\sigma\big(r,u^m(r)\big)\big\|^2_{\text{HS}}\mathrm{d}r\Big)^{\frac{4p}{p-2}}\Big]\\&\leq C_pE\Big[
	\Big(\int_{s}^{t}\Big\|\sigma\big(r,u^m(r)\big)\Big\|^p_{\text{HS}}\mathrm{d}r\Big)^\frac{8}{p-2}|t-s|^4\Big]\\&\leq C_{p,K,T}|t-s|^4.
\end{align*}
By Corollary 1.2 in \cite{W}, there exists a constant $C_3$ depending on ${p,K,T}$ (independent of $k,m$) and a random variable $Y^m_k$  such that there exists  $0<\delta_1<p-2$ such that
\begin{align}\label{I1}
\Big|III_k^m(t,s)\Big|\leq Y^m_k|t-s|^\frac{p-2-\delta_1}{4p},\text{$\forall\ 0\leq s\leq t\leq T$  a.s. }
\end{align}
Moreover,
\begin{align}\label{C3}
E\Big[\big|Y^m_k\big|^\frac{8p}{p-2}\Big]<C_3.
\end{align}
For the term $I_k^m(t,s)$, we have
\begin{eqnarray}\label{II1}
 I_k^m(t,s)&\leq&|t-s|^\frac{1}{2}\Big(\int_{0}^{T}\big|\sum_{j=1}^{m}A_j^m(r)b_{jk}(r)\big|^2\mathrm{d}r\Big)^\frac{1}{2}\nonumber\\
 &\leq&|t-s|^\frac{1}{2}\Big(\int_{0}^{T}\sum_{j=1}^{m}\big|A_j^m(r)\big|^2\cdot\sum_{j=1}^{m}\big|b_{jk}(r)\big|^2\mathrm{d}r\Big)^\frac{1}{2}\nonumber\\
 &\leq& C_k|t-s|^\frac{1}{2}\Big(\int_{0}^{T}\big|u^m(r)\big|^2_r \mathrm{d}r\Big)^\frac{1}{2},
\end{eqnarray}
where the last inequality follows from the fact that
$$\sum_{j=1}^{m}|b_{jk}(r)|^2\leq \sum_{j=1}^{\infty}|b_{jk}(r)|^2$$
is uniformly bounded (with respect to $m$) on the interval $[0,T]$.
The term $II_k^m(t,s)$ can be bounded as
\begin{align}\label{III1}
 II_k^m(t,s)&\leq|t-s|^\frac{1}{2}\Big(\int_{s}^{t}\lambda_k^2(r)\big|A^m_k(r)\big|^2\mathrm{d}r\Big)^\frac{1}{2}\leq|t-s|^\frac{1}{2}C_k\Big(\int_{s}^{t}\big|u^m(r)\big|^2_r \mathrm{d}r\Big)^\frac{1}{2}.
\end{align}
Combining (\ref{I1})-(\ref{III1}), we can find  some constants $C_k$, independent of $m$, such that
 $$\big|A_k^m(t)-A_k^m(s)\big|\leq C_k|t-s|^\frac{1}{2}\Big(\int_{0}^{T}\big|u^m(r)\big|^2_r \mathrm{d}r\Big)^\frac{1}{2}+Y_k^m|t-s|^\frac{p-2-\delta_1}{4p}.$$
For $N>0$ and $q_k>0$, define
\begin{align*}
	J_k^N:=\left\{g\in C([0,T];\mathbb{R}):\big|g(t)-g(s)\big|\leq C_kTN|t-s|^\frac{1}{2}+q_k|t-s|^\frac{p-2-\delta_1}{4p}\right\}.
\end{align*}
Set $K_{N,J}=\cap_{i=1}^\infty K_{N,J_i^N}$.
By Corollary \ref{probest}, for any $\varepsilon>0$, there exists $N>0$ such that $$\sup_{n\in\mathbb{N}}P\Big(\sup_{t\in[0,T]}|u^n(t)|_t>N\Big)\leq\frac{\varepsilon}{4},$$
$$\sup_{n\in\mathbb{N}}P\Big(\int_{0}^{T}\|u^n(t)\|^2_t\mathrm{d}t>N\Big)\leq\frac{\varepsilon}{4}.$$
For the relative compact subset $K_{N,J}$, we have
\begin{eqnarray*}
	&&P\big(u^m\notin K_{N,J}\big)\\
	&\leq& P\Big(\sup_{t\in[0,T]}|u^m(t)|_t>N\Big)+P\Big(\int_{0}^{T}\|u^m(t)\|^2_t\mathrm{d}t>N\Big)+\sum_{k=1}^{\infty}P\big(Y_k^m>q_k\big)\\
	&\leq&\frac{\varepsilon}{2}
+\sum_{i=1}^{\infty}\frac{E\big[(Y_i^m)^\frac{8p}{p-2}\big]}{q_i^\frac{8p}{p-2}} \\
 &\leq&\frac{\varepsilon}{2}
+C_3\sum_{i=1}^{\infty}\frac{1}{q_i^{\frac{8p}{p-2}}},
\end{eqnarray*}
where the last inequality follows from (\ref{C3}). Choose $q_i, i\geq 1$ sufficiently large to obtain
 	$$P\big(u^m\notin K_{N,J}\big)<\varepsilon,\ \forall\ m>0.$$
Since $\varepsilon$ is arbitrary, taking into account Lemma $\ref{tight}$,  we have shown that $\{\mathcal{L}(u^m)\}_{m\geq1}$ is tight as a family of probability measures on $\mathbb{Z}$.$\hfill\blacksquare$

\vskip 0.5cm
By a generalized Skorokhod representation theorem (see Theorem C.1 in \cite{ZEP} or Theorem A.1 in \cite{PKR}), there exists a new probability space $(\Omega^*,\mathcal{F}^*,P^*)$, a sequence of $\mathbb{Z}$-valued random processes $\{u_*,u^m_*,m\geq1\}$  and an $\mathbb{H}$-cylindrical Brownian motion $W^*$ such that
\begin{eqnarray}\label{eq 2023}
P^*\circ(u^m_*,W^*)^{-1}=P\circ(u^m,W)^{-1}
\end{eqnarray}
and that, taking a subsequence if necessary,  $\lim_{m\rightarrow\infty}u_{*}^{m}= u_*$ in $\mathbb{Z}$, $P^*\text{-}a.s$. By the Fatou lemma, we also have $u_*\in \mathbb{L}^2(\Omega^*;L^\infty([0,T];\mathbb{H})\cap\mathbb{Y})$.
Moreover, we have the following stronger convergence result which will be used later.
	\begin{lemma}\label{cauchy}
		$\{u_*^n\}_{n\geq1}$ is a Cauchy sequence in probability in the space $\mathbb{X}\cap\mathbb{Y}$.
	\end{lemma}
	\textbf{Proof.}
	 Without loss of generality, we assume $m>n$.  To simplify the notations,  we  set
 \begin{eqnarray}\label{eq 2023 01}
A_k^{n,*}(t):=\big(u^n_*(t),e_k(t)\big)_t.
 \end{eqnarray}
  Let $k\leq  n$, we have
	\begin{eqnarray*}
		&&A_k^{m,*}(t)-A_k^{n,*}(t)\nonumber\\
		&=&\sum_{j=1}^{n}\int_{0}^{t} \left(A^{m,*}_j(s)-A^{n,*}_j(s)\right) b_{jk}(s) \mathrm{d}s+\sum_{j=n+1}^{m}\int_{0}^{t} A^{m,*}_j(s) b_{jk}(s) \mathrm{d}s\nonumber\\
		&&+\int_{0}^{t}\lambda_k(s)\big(A^{m,*}_k(s)-A^{n,*}_k(s)\big)\mathrm{d}s\nonumber\\
		&&+\int_{0}^{t}\Big(e_k(s),\big(\sigma\big(s,u_*^m(s)\big)-\sigma\big(s,u_*^n(s)\big)\big)\mathrm{d}W^*_s\Big).
	\end{eqnarray*}
The following relationship holds due to the orthogonality of $\{e_j(t)\}_{j\geq1}$.
	\begin{equation}\label{umneq}
		\big|u_*^{m}(t)-u_*^n(t)\big|^2_t=\sum_{k=1}^{n} \big|A_k^{m,*}(t)-A_k^{n,*}(t)\big|^2+\sum_{k=n+1}^{m}\big|A_k^{m,*}(t)\big|^2,
	\end{equation}
	\begin{equation}\label{umneq2}
	\big\|u_*^{m}(t)-u_*^n(t)\big\|^2_t=-\sum_{k=1}^{n} \lambda_k(t)\big|A_k^{m,*}(t)-A_k^{n,*}(t)\big|^2-\sum_{k=n+1}^{m}\lambda_k(t)\big|A_k^{m,*}(t)\big|^2.
\end{equation}

The proof of this lemma is divided into two steps.

{\bf Step 1}. We will prove that
\begin{align}\label{000}
    \sup_{0\leq t \leq T}\sum_{k=1}^{n} |A_k^{m,*}(t)-A_k^{n,*}(t)|^2-2\sum_{k=1}^n \int_0^T \lambda_k(s)\big(A_k^{m,*}(s)-A_k^{n,*}(s)\big)^2 \mathrm{d}s\rightarrow0
\end{align}
 in probability as $m, n\rightarrow\infty$.
\vskip 0.3cm
For $l\geq 1$, denote a  time-dependent orthogonal projection by
$$P^s_l\big(u(s)\big)=\sum_{i=1}^{l}\big(u(s),e_i(s)\big)_se_i(s).$$
By ${\rm It\hat{o}}$'s formula,
	\begin{align}
		&\ \ \sum_{k=1}^{n} |A_k^{m,*}(t)-A_k^{n,*}(t)|^2-2\sum_{k=1}^n \int_0^t \lambda_k(s)\big(A_k^{m,*}(s)-A_k^{n,*}(s)\big)^2 \mathrm{d}s\nonumber\\
		=&\ 2\sum_{k=1}^n\sum_{j=1}^n \int_0^t \big(A_j^{m,*}(s)-A_j^{n,*}(s)\big)b_{jk}(s)\big(A_k^{m,*}(s)-A_k^{n,*}(s)\big) \mathrm{d}s  \nonumber\\
		&+2\sum_{k=1}^n\sum_{j=n+1}^m \int_0^t A_j^{m,*}(s)b_{jk}(s)\big(A_k^{m,*}(s)-A_k^{n,*}(s)\big) \mathrm{d}s  \nonumber\\&+2\int_{0}^{t}\Bigg(P_n^s\big(u_*^m(s)-u_*^n(s)\big),\Big(\sigma\big(s,u_*^m(s)\big)-\sigma\big(s,u_*^n(s)\big)\Big)\mathrm{d}W^*_s\Bigg)\nonumber
\\&+\int_{0}^{t}\sum_{k=1}^{n}\sum_{j=1}^{\infty}\big|\sigma_j^k\big(s,u_*^m(s)\big)-\sigma_j^k\big(s,u_*^n(s)\big)\big|^2\mathrm{d}s\nonumber\\
=:&\ \text{I}^m_n(t)+\text{II}^m_n(t)+\text{III}^m_n(t)+\text{IV}^m_n(t).
	\end{align}
	Since $b_{jk}(s)=-b_{kj}(s)$,
	\begin{align}\label{I}
		\text{I}^m_n(t)=0.
	\end{align}
For the second term, we have 	
	\begin{align*}
		&\ \ \ \ \Big|\text{II}^m_n(t)\Big|\\&=2\Big|\sum_{k=1}^n\sum_{j=n+1}^m \int_0^t A_j^{m,*}(s)\frac{a_s^{\prime}}{a_s}\cdot \frac{2jk}{j^2-k^2}\big(A_k^{m,*}(s)-A_k^{n,*}(s)\big)(-1)^{j+k} \mathrm{d}s\Big|\\
		&=2 \Big|\int_0^t  \sum_{k=1}^n\sum_{j=n+1}^m \frac{A_j^{m,*}(s) j\pi}{a_s}\cdot\frac{\big(A_k^{m,*}(s)-A_k^{n,*}(s)\big)k\pi}{a_s}\cdot \frac{2a_s a_s^{\prime}}{\pi^2(j^2-k^2)}(-1)^{j+k}\mathrm{d}s\Big|\\
		&\leq 2 \int_0^t  \Big(\!\!-\!\!\!\!\sum_{j=n+1}^m |A_j^{m,*}(s)|^2\lambda_j(s)\Big)^{1/2}\Bigg(\!\sum_{j=n+1}^m\!\!\Big(\!\sum_{k=1}^n\big|A_k^{m,*}(s)-A_k^{n,*}(s)\big|\sqrt{-\lambda_k}\frac{2|a_s^\prime| a_s}{\pi^2(j^2-k^2)}\Big)^2\!\Bigg)^{\frac{1}{2}}\! \!\!\mathrm{d}s \\
		&\leq 2\int_0^t \|u^m_*(s)\|_s \Big(-\sum_{k=1}^n\big|A_k^{m,*}(s)-A_k^{n,*}(s)\big|^2 \lambda_k(s) \Big)^{\frac{1}{2}} \Big(\sum_{j=n+1}^{m}\sum_{k=1}^{n}\frac{1}{(j^2-k^2)^2}\Big)^{\frac{1}{2}} \cdot \frac{2|a_s^\prime| a_s}{\pi^2}\mathrm{d}s\\
		&\leq 2\int_0^t \|u^m_*(s)\|_s \|u^m_*(s)-u^n_*(s)\|_s \Big(\sum_{j=n+1}^{m}\sum_{k=1}^{n}\frac{1}{(j^2-k^2)^2}\Big)^{\frac{1}{2}} \cdot \frac{2|a_s^\prime| a_s}{\pi^2} \mathrm{d}s \\
		&\leq \frac{4L^2}{\pi^2} \left(\sum_{j=n+1}^{m}\sum_{k=1}^{n}\frac{1}{(j^2-k^2)^2}\right)^{\frac{1}{2}}
		\|u_*^m\|_{\mathbb{Y}} \|u_*^m-u_*^n\|_{\mathbb{Y}}.
	\end{align*}
	Taking expectation at both side and using Lemma $\ref{estimate}$, we have
	\begin{align*}
		E^*\Big[\sup_{0\leq t \leq T}\Big|\text{II}^m_n(t)\Big|\Big]&\leq\ C_L\Big(\sum_{j=n+1}^{m}\sum_{k=1}^{n}\frac{1}{(j^2-k^2)^2}\Big)^\frac{1}{2}E^*\Big[\|u_*^m\|_\mathbb{Y}^2\Big]^\frac{1}{2}E^*\Big[\|u_*^m-u_*^n\|_\mathbb{Y}^2\Big]^\frac{1}{2}\\&\leq 2C_1C_{K,L}\Big(\sum_{j=n+1}^{m}\sum_{k=1}^{n}\frac{1}{(j^2-k^2)^2}\Big)^\frac{1}{2},
	\end{align*}
where $E^*$ denotes the expectation with respect to $P^*$.
	Now,  setting $l=n+1-k$, we see that
	\begin{align*}
		\sum_{j=n+1}^{m} \sum_{k=1}^{n} \frac{1}{(j^2-k^2)^2}&=\sum_{j=n+1}^{m}\sum _{k=1}^n \frac{1}{(j+k)^2(j-k)^2}\\
		&\leq \sum_{j=n+1}^{m} \frac{1}{j^2} \sum_{k=1}^{n} \frac{1}{(j-k)^2}\\
		&\leq \sum_{j=n+1}^{m} \frac{1}{j^2} \sum_{l=1}^{n} \frac{1}{l^2} \rightarrow 0 \textit{ as m, n } \rightarrow \infty.
	\end{align*}
	It follows that
	\begin{align*}
		&	E^*\Big[\sup_{0\leq t \leq T}\Big|\text{II}^m_n(t)\Big|\Big]\rightarrow 0 \textit{ as m, n }\rightarrow \infty.
	\end{align*}
In particular,
\begin{align}\label{II}
\sup_{0\leq t \leq T}\Big|\text{II}^m_n(t)\Big|\rightarrow 0 \text{ in probability as $m,n\rightarrow\infty$}.
\end{align}
Now, we want to get an estimate for $\text{III}^m_n(t)$.  Notice that
$\text{III}^m_n(t), t\geq0$ is a continuous local martingale such that
\begin{align*}
	\langle \text{III}^m_n \rangle_t\leq 4\int_{0}^{t}|u_*^m(s)-u_*^n(s)|_s^2\big\|\sigma\big(s,u_*^m(s)\big)-\sigma\big(s,u_*^n(s)\big)\big\|^2_{\text{HS}}\mathrm{d}s.
\end{align*}
For $\delta>0$, define	
 $$\tau^{n,m}_\delta:=\inf\big\{t>0:	\langle \text{III}^m_n \rangle_{t}>\delta^2\big\}\wedge T.$$
Then for any $\varepsilon>0$, by the Cheybshev inequality and BDG inequality,
\begin{align*}
	&\ \ \ \ P^*\Bigg(\sup_{0\leq t \leq T}\Big|\text{III}^m_n(t)\Big|\geq\varepsilon\Bigg)\\
&\leq P^*\Big(\tau_\delta^{n,m}=T,\sup_{0\leq t \leq T}\big|\text{III}^m_n(t)\big|\geq\varepsilon\Big)+P^*\Big(\tau_\delta^{n,m}<T,\sup_{0\leq t \leq T}\big|\text{III}^m_n(t)\big|\geq\varepsilon\Big)\\
&\leq\frac{3}{\varepsilon}E^*\Big[\langle \text{III}^m_n\rangle_{\tau^{n,m}_\delta}^\frac{1}{2}\Big]+P^*\Big(\tau^{n,m}_\delta<T\Big)\\ &\leq\frac{3\delta}{\varepsilon}+ P^*\Big(4\int_{0}^{T}|u_*^m(s)-u_*^n(s)|_s^2\big\|\sigma\big(s,u_*^m(s)\big)-\sigma\big(s,u_*^n(s)\big)\big\|^2_{\text{HS}}\mathrm{d}s>\delta^2\Big).
\end{align*}
   Observe that by Proposition $\ref{estimate}$, for $\forall\ \varepsilon_1>0$, there exists $M>0$ such that for $\forall\ m,\ n\in \mathbb{N}$,
 \begin{align*}
 	P^*\Big(\sup_{0\leq t \leq T}\big|u_*^m(t)\big|_t\vee\sup_{0\leq t \leq T}\big|u_*^n(t)\big|_t>M\Big)\leq\varepsilon_1.
 \end{align*}
Then by Assumption \ref{assump2} (i), we have
\begin{eqnarray}\label{1410}
   && P^*\Big(4\int_{0}^{T}\big|u_*^m(s)-u_*^n(s)\big|_s^2\big\|\sigma\big(s,u_*^m(s)\big)-\sigma\big(s,u_*^n(s)\big)\big\|^2_{\text{HS}}\mathrm{d}s>\delta^2\Big)\nonumber\\
   &\leq& P^*\Big(\int_{0}^{T}\big|u_*^m(s)-u_*^n(s)\big|_s^2\mathrm{d}s>\frac{\delta^2}{16M^2K^2}\Big)+\varepsilon_1.
\end{eqnarray}
Since $u_*^m\rightarrow\ u_*$ a.s. in $\mathbb{Z}$, we deduce that
\begin{eqnarray}
	\lim_{m,n\rightarrow\infty}P^*\Bigg(\sup_{0\leq t \leq T}\Big|\text{III}^m_n(t)\Big|>\varepsilon\Bigg)
	\leq\varepsilon_1+\frac{3\delta}{\varepsilon}\nonumber.
\end{eqnarray}
Letting  $\varepsilon_1$ and $\delta$ tend to zero, we obtain that
\begin{align}\label{III}
\sup_{0\leq t \leq T}\Big|\text{III}^m_n(t)\Big|\rightarrow 0\text{ in probability as $m, n\rightarrow\infty$}.
\end{align}
By a similar estimate, we can also show that
	\begin{align}\label{IV}
	\sup_{0\leq t \leq T}\Big|\text{IV}^m_n(t)\Big|\rightarrow0\text{ a.s. as $m,n\rightarrow\ \infty$.}
	\end{align}
	Combining (\ref{I}), (\ref{II}), (\ref{III}), and (\ref{IV}) together, we obtain that
\begin{align*}
    \sup_{0\leq t \leq T}\sum_{k=1}^{n} |A_k^{m,*}(t)-A_k^{n,*}(t)|^2-2\sum_{k=1}^n \int_0^T \lambda_k(s)\big|A_k^{m,*}(s)-A_k^{n,*}(s)\big|^2 \mathrm{d}s\rightarrow\ 0
\end{align*}
as $m, n\rightarrow\infty$ in probability.\\
{\bf Step 2}. We prove that the remaining terms also tend to zero, i.e.,
\begin{align}\label{001}
	\sup_{0\leq t \leq T}\sum_{k=n}^{m} |A_k^{m,*}(t)|^2-2\sum_{k=n}^m \int_0^T \lambda_k(s)\big|A_k^{m,*}(s)\big|^2 \mathrm{d}s\rightarrow\ 0
\end{align}
as $m,n\rightarrow \infty$ in probability.
\vskip 0.3cm
	We only show that the stochastic integral term involved tends to zero in probability since the other terms can be handled similarly as in step 1. \\
	The stochastic integral term is given by:
	$$ \text{V}^m_n(t)=\int_{0}^{t}\Big(\sum_{i=n}^{m}\big(u^m_*(s),e_i(s)\big)_se_i(s),\sigma\big(s,u_*^m(s)\big)\mathrm{d}W^*_s\Big).$$
Following the same arguments as in the proof of (\ref{III}), for arbitrarily small constants $\delta, \varepsilon_1>0$, and sufficiently big constant $M>0$, we deduce that
	\begin{align*}
		&\ \ \ P^*\Bigg(\sup_{0\leq t \leq T}\Big|\text{V}^m_n(t)\Big|>\varepsilon\Bigg)		\\&\leq\frac{3\delta}{\varepsilon}+\varepsilon_1+P^*\Big(\int_{0}^{T}\big|\sum_{i=n}^{m}\big(u^m_*(s),e_i(s)\big)_se_i(s)\big|_s^2\mathrm{d}s>\frac{\delta^2}{K^2(M+1)^2}\Big),
	\end{align*}
while
\begin{eqnarray*}
	&&P^*\Big(\int_{0}^{T}\big|\sum_{i=n}^{m}\big(u^m_*(s),e_i(s)\big)_se_i(s)\big|_s^2\mathrm{d}s>\frac{\delta^2}{K^2(M+1)^2}\Big)\\
	&\leq& P^*\Big(\int_{0}^{T}|u^m_*(s)-u_*(s)|_s^2\mathrm{d}s>\frac{\delta^2}{4K^2(M+1)^2}\Big)\\
	&&+P^*\Big(\int_{0}^{T}\sum_{i=n}^{\infty}\big(e_i(s),u_*(s)\big)_s^2\mathrm{d}s>\frac{\delta^2}{4K^2(M+1)^2}\Big).
\end{eqnarray*}
The first term tends to 0 as $m\rightarrow\infty$  because $\lim_{m\rightarrow\infty}u_*^m= u_*$  a.s. in $\mathbb{Z}$. The second term also tends to 0 by Chebyshev's inequality, the dominated convergence theorem and the fact that $u_*\in\mathbb{L}^2(\Omega^*;\mathbb{Z})$. Thus we conclude that
\begin{align*}
	\Big\{\sup_{0\leq t \leq T}\sum_{k=n}^{m} |A_k^{m,*}(t)|^2-2\sum_{k=n}^m \int_0^T \lambda_k(s)\big(A_k^{m,*}(s)\big)^2 \mathrm{d}s\Big\}\rightarrow\ 0
\end{align*}
in probability as $m,n\rightarrow \infty$.\\

Combining step 1 and step 2, we see that $\{u_*^n\}_{n\geq1}$ is a Cauchy sequence in probability in the space  $\mathbb{X}\cap \mathbb{Y}$.
Hence, we can assume that $u_*\in \mathbb{X}\cap \mathbb{Y}$ and there exists a subsequence ${n_k}$ such that $\lim_{k\rightarrow\infty}u_*^{n_k}= u_*\ \text{in}\ \mathbb{X}\cap \mathbb{Y}$ a.s. $\hfill\blacksquare$
\section{Proof of the main result}
	After all the necessary preparations in Section 3, we are ready to complete the proof of the main result.  Through the proof of Lemma \ref{cauchy}, we find that $\{A^{n,*}:=(A^{n,*}_1,A^{n,*}_2,...,A^{n,*}_n,0,0...)\}_{n\geq1}$ is a Cauchy sequence  in $C([0,T];l^2)$  in probability, which converges to some $(A^*_1,A^*_2,A^*_3,...)$ in $C([0,T];l^2)$. It can be easily seen that the limit $u_*(t,x)$ obtained in Section 3 admits the following expansion
	\begin{eqnarray}
		u_*(t,x)=\sum_{k=1}^{\infty} A^*_k(t) e_k(t,x),
		\end{eqnarray}
	where the  series  converges  in $\mathbb{X}$. Now we prove (\ref{eq 2023 02}). Then  $u_*$ is indeed a weak  solution in probability sense of the stochastic heat equation (\ref{SHE}). We first have the following result.
	\begin{lemma}\label{2.4}
		The stochastic processes $A^*_k,k\geq 1$  satisfy the following  system of equations: for any $k\geq 1$ and $t\geq 0$, \ a.s.
		\begin{eqnarray}\label{00}
			&& A^*_k(t)-(u_0,e_k(0))_0-\sum_{j=1}^{\infty}\int_{0}^{t}A^*_j(s)b_{jk}(s)\mathrm{d}s\\
&=&\int_{0}^{t}\lambda_k(s)A^*_k(s)\mathrm{d}s\nonumber
   +\int_{0}^{t}\Big(e_k(s),\sigma\big(s,u_*(s)\big)\mathrm{d}W^*_s\Big).
		\end{eqnarray}	
	\end{lemma}
	\begin{remark}\label{exchange}
		In the equation above, we can interchange the infinite sum with the integral $$\int_{0}^{t}\sum_{j=1}^{\infty}A^*_j(s)b_{jk}(s)\mathrm{d}s=\sum_{j=1}^{\infty}\int_{0}^{t}A^*_j(s)b_{jk}(s)\mathrm{d}s,$$   because
		\begin{align*}
			&\ E^*\Big[\int_{0}^{T}\sum_{j=1}^{\infty}|A^*_j(s)||b_{jk}(s)|\mathrm{d}s\Big]\\\leq\ &\sum_{j=1}^{\infty}E^*\Big[\int_{0}^{T}|A^*_j(s)|\frac{|a'_s|}{a_s}\frac{2jk}{|j^2-k^2|}\mathrm{d}s\Big]\\=\ &\sum_{j=1}^{\infty}E^*\Big[\int_{0}^{T}\frac{j|A^*_j(s)|}{a_s}|a'_s|\cdot\frac{2k}{|j^2-k^2|}\mathrm{d}s\Big]\\\leq\ &CE^*\Big[\int_{0}^{T}\big(-\sum_{j=1}^{\infty}|A^*_j(s)|^2\lambda_j(s)\big)^\frac{1}{2}\big(\sum_{j=1}^{\infty}\frac{4k^2}{(j^2-k^2)^2}\big)^\frac{1}{2}\mathrm{d}s\Big]\\\leq\ &Ck\sum_{j=1}^{\infty}\Big(\frac{1}{(j^2-k^2)^2}\Big)^\frac{1}{2}E^*\Big[\int_{0}^{T}\|u_*(s)\|_s\mathrm{d}s\Big]\\<\ &\infty. \end{align*}

	\end{remark}
	\textbf{Proof.} According to the equations satisfied by $\{A_k^{n},k=1,2,...,n; W\}$ (see (\ref{approxeq})) and by (\ref{u^n}) (\ref{eq 2023}) (\ref{eq 2023 01}), for $ n\geq 1, k\leq n,$  we have
	\begin{eqnarray*}
		&&A_k^{n,*}(t)-(u_0,e_k(0))_0-\sum_{j=1}^{n}\int_{0}^{t}A^{n,*}_j(s) b_{jk}(s) \mathrm{d}s\\
		&=&\int_{0}^{t}\lambda_k(s)A^{n,*}_k(s)\mathrm{d}s+\int_0^t\Big(e_k(s),\sigma\big(s,u_*^n(s)\big)\mathrm{d}W^*_s\Big).
	\end{eqnarray*}
	We will let $n\rightarrow \infty$ and show that each of the terms in the above equation converges to the corresponding terms in (\ref{00}). According to Lemma 3.6, without loss of generality, we can assume that ${u_*^n}$ converges to $u_*$ a.s. in $\mathbb{X}\cap \mathbb{Y}$ as $n\rightarrow\infty$.  We only need to pay attention to the convergence of the third term on the left hand side, i.e.
	\begin{align}\label{01}
		\sum_{j=1}^{n}\int_{0}^{\cdot}A_j^{n,*}(s)b_{jk}(s)\mathrm{d}s\rightarrow\int_{0}^{\cdot}\sum_{j=1}^{\infty}A^*_j(s)b_{jk}(s)\mathrm{d}s
	\end{align}
	 in $C([0,T];\mathbb{R})$ in probability as $n\rightarrow\infty$. The convergence of the remaining terms is obvious.

The proof of (\ref{01}) is  divided into two parts,
	\begin{itemize}
		\item [(i)] $$\sum\limits_{j=1}^{n}\int_{0}^{T}|A_j^{n,*}(s)-A_j^*(s)||b_{jk}(s)|\mathrm{d}s\rightarrow\ 0,\text{ a.s. as $n\rightarrow\infty$,}$$
		\item [(ii)] $$E^*\Big[\sum\limits_{j=n+1}^{\infty}\int_{0}^{T}|A_j^*(s)||b_{jk}(s)|\mathrm{d}s\Big]\rightarrow\ 0,\text{ as $n\rightarrow\infty$.}$$
	\end{itemize}
	Taking into account Remark $\ref{exchange}$, both (i) and (ii) can be proved using the similar arguments as in the proof of Lemma \ref{cauchy}. So we omit the details here.$\hfill\blacksquare$
\vskip 0.3cm

   Now, we show that
	\begin{proposition}\label{existence}
		The random field $u_*$ obtained above is a solution to the stochastic heat equation (\ref{SHE}), namely, for all $t\in[0,T]$ and $\varphi\in\mathcal{C}_0^\infty (\bar{\mathcal{O}}_T)$,
		\begin{align}\label{eq 2023 02}
			&\int_{0}^{a_t} u_*(t,x)\varphi(t,x)\mathrm{d}x -\int_{0}^{a_0} u_0(x)\varphi(0,x) \mathrm{d}x-\int_{0}^{t}\int_{0}^{a_s} u_*(s,x)\varphi'(s,x) \mathrm{d}x\mathrm{d}s\nonumber\\
			=&\int_{0}^{t}\int_{0}^{a_s}  u_*(s,x)\Delta\varphi(s,x)\mathrm{d}x\mathrm{d}s+\int_0^t\Big(\varphi (s),\sigma\big(s,u_*(s)\big)dW^*_s\Big)\ a.s.
		\end{align}
	\end{proposition}
	\textbf{Proof.}
 The statement of the first sentence of this section shows that the proof of this proposition is complete once we have proved (\ref{eq 2023 02}).
In the rest of the proof we shall verify  (\ref{eq 2023 02}).

	Let $\varphi$ be a test function, i.e. $\varphi\in\mathcal{C}_0^\infty (\bar{\mathcal{O}}_T)$  and define $\varphi_k(s):=\big(\varphi(s),e_k(s)\big)_s$. Since $e_k(s,x)$ and $\varphi(s,x)$ vanish at the boundary of $I_s$, we have
	\begin{align*}
		\mathrm{d}\varphi_k(s)=\big(\varphi'(s),e_k(s)\big)_s\mathrm{d}s+\big(\varphi(s),e_k'(s)\big)_s\mathrm{d}s.
	\end{align*}
	Recall that $\mathrm{d}A^*_k(s)=\big(u_*(s),e_k'(s)\big)_s\mathrm{d}s+\lambda_k(s)A^*_k(s)\mathrm{d}s+\Big(e_k(s),\sigma\big(s,u_*(s)\big)\mathrm{d}W^*_s\Big).$
	By the ${\rm It\hat{o}}$ formula, we have
	\begin{align*}
		\mathrm{d}\varphi_kA^*_k(s)=&\ \varphi_k(s)\mathrm{d}A_k^*(s)+A_k^*(s)\mathrm{d}\varphi_k(s)\\=&\ \varphi_k(s)\Bigg(\big(u_*(s),e_k'(s)\big)_s\mathrm{d}s+\lambda_k(s)A_k^*(s)\mathrm{d}s+\Big(e_k(s),\sigma\big(s,u_*(s)\big)\mathrm{d}W^*_s\Big)\Bigg)\\&\ +A_k^*(s)\Big(\big(\varphi'(s),e_k(s)\big)_s+\big(\varphi(s),e_k'(s)\big)_s\Big)\mathrm{d}s,
	\end{align*}
	that is,
	\begin{align}\label{to sum} &\ \varphi_kA^*_k(t)\nonumber\\=&\  \varphi_k(u_0,e_k(0))_0+\int_{0}^{t}\varphi_k(s)\big(u_*(s),e'_k(s)\big)_s\mathrm{d}s+\int_{0}^{t}\lambda_k(s)\varphi_k(s)A^*_k(s)\mathrm{d}s\nonumber\\
&\  +\int_{0}^{t}\Big(\varphi_k(s)e_k(s),\sigma\big(s,u_*(s)\big)\mathrm{d}W^*_s\Big)+\int_{0}^{t}A^*_k(s)\big(\varphi'(s),e_k(s)\big)_s\mathrm{d}s\nonumber\\
     &\ +\int_{0}^{t}A^*_k(s)\big(\varphi(s),e_k'(s)\big)_s\mathrm{d}s\nonumber\\
     =&\ \mathbb{I}_k+\mathbb{II}_k(t)+\mathbb{III}_k(t)+\mathbb{IV}_k(t)+\mathbb{V}_k(t)+\mathbb{VI}_k(t).
	\end{align}
	Recall that $u_*(t,x)=\sum\limits_{j=1}^{\infty}A^*_j(t)e_j(t,x)$ and  $\big(u_*(t),\varphi(t)\big)_t=\sum\limits_{k=1}^{\infty}\varphi_kA^*_k(t)$.
	Adding up (\ref{to sum}) to an arbitrarily big natural number $n$ and then letting $n\rightarrow \infty$, we obtain
	\begin{itemize}
		\item [(i)] $$\sum\limits_{k=1}^{\infty}\mathbb{I}_k=\big(\varphi(0),u(0)\big)_0,$$
		\item[(ii)]
		$$\sum\limits_{k=1}^{\infty}\mathbb{III}_k(t)=\int_{0}^{t}\Delta\varphi(s)u_*(s)\mathrm{d}s\ a.s.,$$
		\item[(iii)]
		$$\sum\limits_{k=1}^{\infty}\mathbb{IV}_k(t)=\int_{0}^{t}\Big(\varphi(s),\sigma\big(s,u_*(s)\big)\mathrm{d}W^*_s\Big) \text{ in probability},$$
		\item[(iv)]
		$$\sum\limits_{k=1}^{\infty}\mathbb{V}_k(t)=\int_{0}^{t}\big(\varphi'(s),u_*(s)\big)_s\mathrm{d}s \text{ in } \mathbb{L}^2(\Omega^*).$$
	\end{itemize}
Here we set $\sum\limits_{k=1}^{\infty}X_k:=\lim_{n\rightarrow\infty}\sum\limits_{k=1}^{n}X_k.$

On the other hand, we denote that
	\begin{align*}
		\mathbb{II}_k(t)=\int_{0}^{t}\varphi_k(s)\big(u_*(s),e'_k(s)\big)_s\mathrm{d}s=\int_{0}^{t}\varphi_k(s)\sum\limits_{j=1}^{\infty}A^*_j(s)b_{jk}(s)\mathrm{d}s.
	\end{align*}
\vskip -0.5cm
	\begin{align*}
         \mathbb{VI}_k(t)=\int_{0}^{t}A^*_k(s)\big(\varphi(s),e'_k(s)\big)_s\mathrm{d}s=\int_{0}^{t}A^*_k(s)\sum\limits_{j=1}^{\infty}\varphi_j(s)b_{jk}(s)\mathrm{d}s.
	\end{align*}
Since $\{b_{jk}(s)\}_{j,k}$ is  skew-symmetric with respect to  $(j,k)$, we see that
$$\sum\limits_{k=1}^{\infty}\mathbb{II}_k(t)+\sum\limits_{k=1}^{\infty}\mathbb{VI}_k(t)=0.$$
To obtain all of the above equations, we have used the interchange of the infinite sum with the integral, which can be seen as follows
	\begin{align}\label{01080247}
		&E^*\Big[\int_{0}^{T}\sum_{k,j}\big|\varphi_k(s)\big|\big|A^*_j(s)\big|\big|b_{jk}(s)\big|\mathrm{d}s\Big]\\
		=\ &\nonumber E^*\Big[\int_{0}^{t}\sum_{k\neq j}\big|\varphi_k(s)\big|\frac{k}{|a_s|}\big|A^*_j(s)\big|\frac{j}{|a_s|}\frac{|a_s||a'_s|}{|j^2-k^2|}\mathrm{d}s\Big]\\
		\leq\ &\nonumber CE^*\Big[\int_{0}^{t}\!\Big(\!\sum_{k}\big|\varphi_k(s)\big|^2\frac{k^2}{|a_s|^2}\sum_{j\neq k}\frac{1}{|j^2-k^2|}\!\Big)^{\frac{1}{2}}\Big(\!\sum_{k}\big|A^*_j(s)\big|^2\frac{j^2}{|a_s|^2}\sum_{j\neq k}\frac{1}{|j^2-k^2|}\!\Big)^\frac{1}{2}\!\mathrm{d}s\Big]\\
		\leq\ &\nonumber CE^*\Big[\int_{0}^{t}\big\|\varphi(s)\big\|_s^2\mathrm{d}s\Big]^\frac{1}{2}E^*\Big[\int_{0}^{t}\big\|u_*(s)\big\|_s^2\mathrm{d}s\Big]^\frac{1}{2}
		\\<\ &\nonumber\infty.
	\end{align}
Putting the above equations together we finally arrive at
 \begin{align}
			&\int_{0}^{a_t} u_*(t,x)\varphi(t,x)\mathrm{d}x -\int_{0}^{a_0} u_0(x)\varphi(0,x) \mathrm{d}x-\int_{0}^{t}\int_{0}^{a_s} u_*(s,x)\varphi'(s,x) \mathrm{d}x\mathrm{d}s\nonumber\\
			=&\int_{0}^{t}\int_{0}^{a_s}  u_*(s,x)\Delta_s\varphi(s,x)\mathrm{d}x\mathrm{d}s+\int_0^t\Big(\varphi (s),\sigma\big(s,u_*(s)\big)\mathrm{d}W^*_s\Big)_{\mathbb{H}}  \ a.s.,
		\end{align}
completing the proof of (\ref{eq 2023 02}).

The proof of this proposition is complete.
	$\hfill\blacksquare$
	\vskip 0.3cm
		Next result is an energy identity/$\rm{It\hat{o}}$-type formula for the solution.
	\begin{proposition}\label{Ito}
		Let $(\tilde{\Omega},\tilde{\mathcal{F}},\{\tilde{\mathcal{F}_t}\}_{t\geq0},\tilde{P},\tilde{W},\tilde{u})$ be a solution  to equation (\ref{SHE}), that is, $(\tilde{\Omega},\tilde{\mathcal{F}},\{\tilde{\mathcal{F}_t}\}_{t\geq0},\tilde{P})$ is an filtered probability space satisfying the usual conditions,  $\tilde{W}$ is an $\mathbb{H}$-cylindrical Wiener process on the given probability space, and  $\tilde{u}$ is a solution of equation (\ref{SHE}) in the sense of Definition \ref{sol1} with $W$ replaced  by $\tilde{W}$. We have
		\begin{align*}
			|\tilde{u}(t)|^2_t=&\ |u(0)|_0^2-2\int_{0}^{t}\|\tilde{u}(s)\|_s^2\mathrm{d}s+2\int_{0}^{t}\Big(\tilde{u}(s),\sigma\big(s,\tilde{u}(s)\big)\mathrm{d}\tilde{W}_s\Big)\\
			&+\int_{0}^{t}\sum_{k=1}^{\infty}\sum_{j=1}^{\infty}\big|\sigma_j^k\big(s,\tilde{u}(s)\big)\big|^2\mathrm{d}s
			\\=&\ |u(0)|_0^2-2\int_{0}^{t}\|\tilde{u}(s)\|_s^2\mathrm{d}s+2\int_{0}^{t}\Big(\tilde{u}(s),\sigma\big(s,\tilde{u}(s)\big)\mathrm{d}\tilde{W}_s\Big)\\&+\int_{0}^{t}\big\|\sigma\big(s,\tilde{u}(s)\big)\big\|_{\rm{HS}}^2\mathrm{d}s.
		\end{align*}
	\end{proposition}	
	\textbf{Proof.}
	Take $\varphi=e_k$ in Definition \ref{sol1}, and define $\tilde{A}_k(t):=\big(\tilde{u}(t),e_k(t)\big)_t$ we have
\begin{align*}
	&\ \ \ \ \tilde{A}_k(t)-{A}_k(0)-\int_{0}^{t}\int_{0}^{a_s} \tilde{u}(s,x) e_k'(s,x) \mathrm{d}x\mathrm{d}s
	\\&= \int_{0}^{t}\lambda_k(s)\tilde{A}_k(s)\mathrm{d}s+\int_{0}^{t}\Big(e_k(s),\sigma\big(s,\tilde{u}(s)\big)\mathrm{d}\tilde{W}_s\Big)_\mathbb{H}.
\end{align*}
By ${\rm It\hat{o}}$'s formula and using similar arguments to prove (\ref{01080247}),
\begin{align*}
	|\tilde{A}_k(t)|^2 =&|(u_0,e_k(0))_0|^2
	+\sum_{j=1}^{\infty}\int_{0}^{t}\tilde{A}_j(s) b_{jk}(s)\tilde{A}_k(s) \mathrm{d}s+2 \int_0^t \lambda_k(s)|\tilde{A}_k(s)|^2 \mathrm{d}s\ \\ &+ 2\int_{0}^{t}\Big(\tilde{A}_k(s)e_k(s),\sigma\big(s,\tilde{u}(s)\big)\mathrm{d}\tilde{W}_s\Big)+\int_{0}^{t}\sum_{j=1}^{\infty}\big|\sigma_j^k\big(s,\tilde{u}(s)\big)\big|^2\mathrm{d}s.
\end{align*}
By adding up $k$ from $1$ to $\infty$ and following the same proof of Lemma $\ref{cauchy}$, we obtain the desired identity.
$\hfill\blacksquare$\par
	\vskip 0.5cm
The next result is about the uniqueness of the solution.
	\begin{theorem}\label{uniqueness}
		The solution of the stochastic heat equation (\ref{SHE}) is pathwise  unique in the space $\mathbb{X}\cap \mathbb{Y}$.
	\end{theorem}
	\textbf{Proof.}
Assume that $u$ and $v$ are two  solutions in $\mathbb{X}\cap \mathbb{Y}$ with initial value $u_0$. Let $w(t):=u(t)-v(t)$. By the energy identity stated in  Proposition \ref{Ito}, we have
	\begin{align}\label{02}
		|w|_t^2=&-2\int_{0}^{t}\|w\|_s^2\mathrm{d}s+2\int_{0}^{t}\Bigg(w(s),\Big(\sigma\big(s,u(s)\big)-\sigma\big(s,v(s)\big)\Big)\mathrm{d}W_s\Bigg)\nonumber\\
&+\int_{0}^{t}\big\|\sigma\big(s,u(s)\big)-\sigma\big(s,v(s)\big)\big\|_{\mathrm{HS}}^2\mathrm{d}s.
	\end{align}
For $M>0$, define
	$$\tau_M:=\inf\Big\{t\in(0,T]:|u(t)|_t\vee|v(t)|_t>M\Big\}\wedge T.$$
Replace $t$ by $t\wedge\tau_M$ and take expectation to get
   	\begin{align*}
   	&\ E\ \big|w(t\wedge\tau_M)\big|_{t\wedge\tau_M}^2\\\leq&-2E\int_{0}^{t\wedge\tau_M}\|w(s)\|_s^2\mathrm{d}s
   +E\int_{0}^{t\wedge\tau_M}\big\|\sigma\big(s,u(s)\big)-\sigma\big(s,v(s)\big)\big\|_{\mathrm{HS}}^2\mathrm{d}s\\   	\leq&-2E\int_{0}^{t\wedge\tau_M}\|w(s)\|_s^2\mathrm{d}s+C_KE\int_{0}^{t\wedge\tau_M}|w(s)|^2_s\mathrm{d}s.
   \end{align*}
	Since $w(\cdot)\in \mathbb{X}\cap \mathbb{Y}$, by the Gronwall inequality, we have $w(\cdot\wedge\tau_M)=0$ a.s. in $\mathbb{X}\cap \mathbb{Y}$. Letting $M\rightarrow\infty$, we have $w=0$ a.s. in $\mathbb{X}\cap \mathbb{Y}$. $\hfill\blacksquare$

	\vskip 0.25cm
\noindent {\bf Completion of the proof of Theorem \ref{WP}}\par
	Proposition \ref{existence} gives the existence of a probabilistic weak solution. Now  Theorem\ \ref{WP} follows from the pathwise uniqueness proved in Theorem \ref{uniqueness} and the  well-known Yamada-Watanabe theorem.$\hfill\blacksquare$
	\vskip 0.5cm
	Given the existence of a unique solution $u$ of equation (\ref{SHE}), the next result shows that the solution $u$ can be approximated by solutions $u^n(t)$ of the finite dimensional systems (\ref{approxeq}) or (\ref{u^n}).

\begin{proposition}\label{Prop 1}
	The solutions $\{u^n\}_{n\geq1}$ of the finite-dimensional interacting systems in (\ref{approxeq}) converge to $u$ in $L^2(\Omega;\mathbb{X}\cap \mathbb{Y})$ as $n\rightarrow\infty$.
\end{proposition}
\textbf{Proof.} The proof is a modification of the proof of Lemma \ref{cauchy}. Recall
$$A_k^n(t)=\big(u^n(t),e_k(t)\big)_t,\quad\quad\quad A_k(t)=\big(u(t),e_k(t)\big)_t.$$
Then, we have
 \begin{equation}\label{03}
 |u(t)-u^n(t)|^2_t=\sum_{k=1}^{n}|A_k(t)-A_k^n(t)|^2+\sum_{k=n+1}^{\infty}|A_k(t)|^2,
 \end{equation}
 \begin{equation}\label{04}
 \|u(t)-u^n(t)\|^2_t=-\sum_{k=1}^{n}\lambda_k(t)|A_k(t)-A_k^n(t)|^2-\sum_{k=n+1}^{\infty}\lambda_k(t)|A_k(t)|^2.
 \end{equation}
 Now by the equations satisfied by $A_k(t), A_k^n(t)$ and Ito's formula,
 \begin{align}\label{05}
 	&\sum_{k=1}^{n}\big|A_k(t)-A_k^n(t)\big|^2-2\sum_{k=1}^{n}\int_{0}^{t}\lambda_k(s)\big|A_k(s)-A_k^n(s)\big|^2\mathrm{d}s\nonumber\\=&\  2\sum_{k=1}^{n}\sum_{j=1}^{n}\int_{0}^{t}\big(A_j(s)-A_j^n(s)\big)b_{jk}(s)\big(A_k(s)-A_k^n(s)\big)\mathrm{d}s\nonumber\\&+\ 2\sum_{k=1}^{n}\sum_{j=n+1}^{\infty}\int_{0}^{t}A_j(s)b_{jk}(s)\big(A_k(s)-A_k^n(s)\big)\mathrm{d}s\nonumber\\&+\ 2\int_{0}^{t}\Bigg(P_n^s\big(u(s)-u^n(s)\big),\Big(\sigma\big(s,u(s)\big)-\sigma\big(s,u^n(s)\big)\Big)\mathrm{d}W_s\Bigg)
 	\nonumber\\&+\ \int_{0}^{t}\sum_{k=1}^{n}\sum_{j=1}^{\infty}\big|\sigma^k_j\big(s,u(s)\big)-\sigma^k_j\big(s,u^n(s)\big)\big|^2\mathrm{d}s\nonumber\\=&\ \text{I}^n(t)+\text{II}^n(t)+\text{III}^n(t)+\text{IV}^n(t).
 \end{align}
 The first term $\text{I}^n(t)=0$  since ${b_{jk}(s)}$ is skew-symmetric with respect to $(j,k)$. Arguing  as in the proof of  Lemma \ref{cauchy}, we have
 \begin{align*}
 	\Big|\text{II}^n(t)\Big|\leq\frac{4L^2}{\pi^2}\Big(\sum_{j=n+1}^{\infty}\sum_{k=1}^{n}\frac{1}{(j^2-k^2)^2}\Big)^\frac{1}{2}\|u\|_\mathbb{Y}\|u-u^n\|_\mathbb{Y}.
 \end{align*}
 This implies, as in the proof of  Lemma \ref{cauchy}, that
 \begin{equation}\label{IIt}
 	E\big[\sup_{0\leq t \leq T}\Big|\text{II}^n(t)\Big|\big]=o(1),\text{ as $n\rightarrow\infty$}.
 \end{equation}
By the BDG inequality, Young inequality and Assumption \ref{assump2} (i),
 \begin{align}\label{IIIt}
 	E\Big[\sup_{0\leq s \leq t}|\text{III}^n(t)|\Big]&\leq C E\Big[\Big(\int_{0}^{t}|u(s)-u^n(s)|_s^2\big\|\sigma\big(s,u(s)\big)-\sigma\big(s,u^n(s)\big)\big\|^2_{\text{HS}}\mathrm{d}s\Big)^\frac{1}{2}\Big]\nonumber\\&\leq C_KE\Big[\sup_{0\leq s \leq t}|u(s)-u^n(s)|_s\Big(\int_{0}^{t}|u(s)-u^n(s)|_s^2\mathrm{d}s\Big)^\frac{1}{2}\Big]\nonumber\\&\leq \frac{1}{2}E\Big[\sup_{0\leq s \leq t}|u(s)-u^n(s)|_s^2\Big]+C_KE\Big[\int_{0}^{t}|u(s)-u^n(s)|_s^2\mathrm{d}s\Big].
 \end{align}
 By Assumption \ref{assump2} (i), we have
  \begin{align}\label{IVt}
 	E\Big[\sup_{0\leq t \leq T}\text{IV}^n(t)\Big]\leq C_K\int_{0}^{T}E\big[|u(s)-u^n(s)|_s^2\big]\mathrm{d}s.
 \end{align}
On the other hand, we also claim that as $n\rightarrow\infty$,
\begin{equation}\label{06}
	E\Big[\sup_{0\leq t \leq T}\sum_{k=n+1}^{\infty}|A_k(t)|^2\Big]\rightarrow 0,
\end{equation}
\begin{equation}
E\Big[\sum_{k=n+1}^{\infty}\int_{0}^{T}\lambda_k(t)|A_k(t)|^2\mathrm{d}t\Big]\rightarrow 0.\label{es}
\end{equation}
(\ref{es}) is easily seen by the dominated convergence theorem and the fact that $u\in L^2(\Omega;\mathbb{Y})$. To prove (\ref{06}), we notice that for a.s. $\omega\in\Omega$ fixed, $$\Big\{ \sum_{k=n+1}^{\infty}|A_k(t)|^2\Big\}_{n\geq 1}$$ is a decreasing sequence of continuous functions on $[0, T]$, which tends to zero for every $t\in [0,T]$. This implies  that as $n\rightarrow\infty$ a.s  $$\sup_{0\leq t \leq T}\sum_{k=n+1}^{\infty}|A_k(t)|^2\rightarrow0.$$
($\ref{es}$) follows now from the dominated convergence theorem.
Combining (\ref{03})-(\ref{es}) together, we arrive at
$$E\Big[\sup_{0\leq s \leq t}|u^n(s)-u(s)|_s^2+\int_{0}^{t}\|u^n(s)-u(s)\|^2_s\mathrm{d}s\Big]\leq\ o(1)\ +\ C_K\int_{0}^{t}E\big[|u^n(s)-u(s)|_s^2\big]\mathrm{d}s.$$
By Gronwall's inequality, we have as $n\rightarrow\infty$,
$$E\Big[\sup_{0\leq t \leq T}|u^n(t)-u(t)|_t^2+\int_{0}^{T}\|u^n(s)-u(s)\|^2_s\mathrm{d}s\Big]\rightarrow 0. $$

The proof of Proposition \ref{Prop 1} is complete.

\vskip 0.4cm
\noindent{\bf Acknowledgement}   This work is partially supported by National Key R\&D Program of China(No. 2022YFA1006001), National Natural Science Foundation of China (Nos. 12131019, 11971456, 11721101). Jianliang Zhai’s research is also supported by the School Start-up Fund (USTC) KY0010000036 and the Fundamental Research Funds for the Central Universities (No. WK3470000016).
	
\end{document}